\pdfoutput=1
\documentclass[11pt]{article}
\usepackage{fullpage}
\usepackage{amsmath, amssymb, amsthm}
\usepackage[utf8]{inputenc}
\usepackage[english]{babel}
\usepackage[numbers]{natbib}
\usepackage{url}
\usepackage[usenames]{color}
\usepackage[colorlinks=true]{hyperref}
\usepackage{cleveref}
\usepackage{tabularx}
\theoremstyle{plain}

\newtheorem{conjecture}{Conjecture}[section]

\newtheorem{lemma}{Lemma}[section]

\newtheorem{theorem}{Theorem}[section]
\newtheorem*{theorem*}{Theorem}

\crefname{conjecture}{Conjecture}{Conjectures}
\crefname{theorem}{Theorem}{Theorems}
\crefname{theorem*}{Theorem}{Theorems}
\crefname{corollary}{Corollary}{Corollaries}
\crefname{lemma}{Lemma}{Lemmas}
\crefname{proposition}{Proposition}{Propositions}
\crefname{remark}{Remark}{Remarks}
\crefname{note}{Note}{Notes}
\crefname{definition}{Definition}{Definitions}
\crefname{notation}{Notation}{Notations}
\crefname{example}{Example}{Examples}
\crefname{question}{Question}{Questions}
\crefname{section}{\S}{Sections}
\newcommand{\floor}[1]{\left\lfloor #1 \right\rfloor}
\newcommand{\Z}{\mathbb{Z}}

\newcommand{\eval}[2]{\left . #1 \right|_{#2}}
\newcommand{\seqnum}[1]{\href{https://oeis.org/#1}{\rm \underline{#1}}}

\begin{document}

\title{Polynomial quotient rings and Kronecker substitution for deriving combinatorial identities}
\author{Joseph M. Shunia \footnote{Post-baccalaureate Student, Department of Mathematics, The Johns Hopkins University, Baltimore, MD, USA;} \footnote{Veeam Software, Columbus, OH, USA (Remote). E-mail: {\tt jshunia1@jh.edu}, {\tt jshunia@gmail.com}}}
\date{March 2024 \\ \small Revised: November 2024, Version 6 \normalsize}

\maketitle

\begin{abstract} \noindent
We introduce a new approach for generating combinatorial identities and formulas by the application of Kronecker substitution to polynomial expansions within quotient rings. Our main result enables the derivation of elementary arithmetic formulas for many C-recursive integer sequences directly from their characteristic polynomials. As sample applications, we present new formulas for the Pell numbers and central binomial coefficients, which are famous integer sequences. These applications lead us to the discovery of a new and unusual formula for the real $n$-th roots of positive integers, $\sqrt[n]{a}$, characterized as the limit of a quotient involving modular exponentiations. From this limit formula we conjecture a fixed-length elementary closed form expression for $\floor{\sqrt[n]{a}}$.
 \\[2mm]
\end{abstract}

\section{Introduction}
Kronecker substitution, named after the mathematician Leopold Kronecker, is a technique that allows for the efficient multiplication polynomials by encoding them as integers in a large base \cite{gathen2013modern}. While this technique has been widely used in the design of fast multiplication algorithms \cite{harvey2009kronecker, harvey2019faster, albrecht2018implementing, bos2020postquantum, greuet2022modular}, its potential applications in combinatorics and number theory have remained largely unexplored until recently \cite{shunia2023simple}.

A \textbf{C-recursive integer sequence of order $d$} is an integer-valued sequence satisfying a recurrence relation with constant coefficients of the form:
\begin{align*}
A(n) = c_{d-1} A(n-1) + c_{d-2} A(n-2) + \cdots + c_{0} A(n-d) ,
\end{align*}
for all $n \geq 0$, where the initial starting conditions $A(0),A(1),\ldots,A(d-1)$ are specified as integers.

In this work, we develop a methodology for applying Kronecker substitution to polynomial expansions within quotient rings, which is useful for generating combinatorial formulas and identities defined by C-recursive integer sequences. \cref{proof:crecursiveencoding} establishes a connection between the coefficients of a polynomial remainder and the integers obtained by evaluating the polynomials at specific values. By defining the modulus to correspond with the characteristic polynomial of C-recursive integer sequences, we can generate new identities for combinatorial sequences using exponential polynomials and Kronecker substitution.

As some sample applications, we prove new formulas for the Pell sequence in \cref{section:pell} and the central binomial coefficients in \cref{section:cbc}. We also deduce in \cref{proof:roots} a remarkable new formula for the real $n$-th roots of positive integers $\sqrt[n]{a}$, which is characterized as the limit of a quotient involving modular exponentiations. From this expression, we conjecture a fixed-length elementary closed-form expression for $\floor{\sqrt[n]{a}}$ (\cref{conjecture:integerroots}).

\section{Preliminaries} \label{section:encodingtheorem}
Given a polynomial $f(x) = a_d x^d + a_{d-1} x^{d-1} + \cdots + a_0$, we denote by $\tilde{f(x)}$ the \textbf{polynomial normalized form} of $f(x)$, given by:
\begin{align*}
\tilde{f}(x) = \frac{f(x)}{a_d} = x^d + \frac{a_{d-1}}{a_d} x^{d-1} + \cdots + \frac{a_0}{a_d} .
\end{align*}

The \textbf{polynomial encoding theorem} was proved by Shunia in \cite{shunia2023simple}. We restate it here as a lemma, as it will be needed for our results:

\begin{lemma} (Shunia, \citep[Theorem 2.1]{shunia2023simple}) \label{proof:encoding}
Given three non-negative integers $b > 0$, $k$, $r \geq k$, and a non-constant polynomial $f(x)$ of non-negative integer coefficients, degree $r$, such that $f(b) \not= 0$, we have that
\begin{align*}
[x^k]f(x) = \floor{\frac{f(f(b))}{f(b)^{k}}} \bmod{f(b)} .
\end{align*}
\end{lemma}
\begin{proof}
The proof is due to Shunia \cite{shunia2023simple} .
\end{proof}

\section{Main results}
We begin with a lemma:

\begin{lemma} \label{proof:recurrences}
Let $n,d \in \Z^+$. Define $A(n)$ to be a C-recursive integer sequence of order $d$ with initial starting conditions $A(0)=A(1)=\cdots=A(d-1)=1$, such that:
\begin{align*}
    A(n) = c_{d-1} A(n-1) + c_{d-2} A(n-2) + \cdots + c_{0} A(n-d) ,
\end{align*}
where the $c_j$ are coefficients in $\Z$. Let $g(x)$ be the characteristic polynomial of $A(n)$, which has the form:
\begin{align*}
    g(x) := c_{d-1} x^{d-1} + \cdots + c_1 x + c_0 \in \Z[x] .
\end{align*}
Fix a ring $R = \Z[x]/I$, where the ideal $I = \langle x^d-g(x) \rangle$. Consider $f(x) := x^n \in R$ to be the image of $x^n$ in $R$. Then $A(n) = f(1)$.
\end{lemma}
\begin{proof}
In the ring $R$, any power of $x^k$ with $k \geq d$ is implicitly modulo the ideal $I = \langle x^d-g(x) \rangle$, which enforces the relation $x^d = g(x)$. For a general $k$, the image of $x^k$ in $R$ can be determined by repeatedly applying the relation
\begin{align*}
x^k = x^{k-1} x = (x^{k-1-d} x^d) x = (x^{k-1-d} g(x)) x = \cdots ,
\end{align*}
until all terms are in terms of $x^{d-1}$ or lower. The image of $x^k$ is thus a linear combination of $\{1, x, x^2, \ldots, x^{d-1}\}$. 

The sequence $A(n)$ is governed by the recurrence relation:
\begin{align*}
A(n) = c_{d-1} A(n-1) + c_{d-2} A(n-2) + \cdots + c_{0} A(n-d) .
\end{align*}
with initial conditions $A(0) = A(1) = \cdots = A(d-1) = 1$ and the corresponding characteristic polynomial $g(x)$ represents these relations
\begin{align*}
g(x) = c_0 + c_1 x + \cdots + c_{d-1} x^{d-1} .
\end{align*}
Under this setup, $x^d = g(x)$ implies that the recurrence relations used in the definition of $R$ match those used to calculate $A(n)$. Thus, the coefficients $a_0, a_1, \ldots, a_{d-1}$ of the image of $f(x)$ in $R$ directly correspond to the coefficients in the linear combination that expresses $A(n)$ in terms of the base cases $A(0), A(1), \ldots, A(d-1)$. This implies that the evaluation at $x=1$ yields the same combination that one would compute directly using the recurrence relation on $A(n)$. Hence, $A(n) = f(1)$.
\end{proof}

\begin{theorem} \label{proof:kroneckerqrings}
Let $k, d \in \Z^+ : k \geq d$. Let $f(x) \in \Z[x]$ such that $f(x)$ is non-constant. Consider a polynomial
\begin{align*}
g(x) := a_d x^d - a_{d-1} x^{d-1} - \cdots - a_0 \in \Z[x] ,
\end{align*}
and the remainder
\begin{align*}
r(x) := f(x)^k \bmod{\tilde{g}(x)} ,
\end{align*}
Let $\gamma, b \in \Z^+$ such that $\gamma^k > |r(b)|$ and $r(b) \not\equiv 0 \pmod{\gamma^k - b}$. Then
\begin{align*}
r(b) = \left( f(\gamma^k)^k \bmod{\tilde{g}(\gamma^k)} \right) \bmod{(\gamma^k - b)} .
\end{align*}
\end{theorem}
\begin{proof}
First, consider the evaluation 
\begin{align*}
    \eval{r(x)}{x=b} = r(b) .
\end{align*}
By the remainder theorem, evaluating a polynomial $h(x) \in \Z[x]$ at $x=b$ is the same as reducing $h(x)$ modulo $(x - b)$. Since we are working modulo $\tilde{g}(x)$, we have the relation
\begin{align*}
    r(b) = \left(f(x)^k \bmod{\tilde{g}(x)}\right) \bmod{(x - b)} .
\end{align*}
Applying Kronecker substitution to all polynomials in the above equation by the substitution $x = \gamma^k$ yields
\begin{align*}
    r(b) = \left(f(\gamma^k)^k \bmod{\tilde{g}(\gamma^k)}\right) \bmod{(\gamma^k - b)} ,
\end{align*}
which is the formula we aimed to prove. Furthermore, recall that we are given $\gamma$ such that $\gamma^k > |r(b)|$. By \cref{proof:encoding}, this implies that the chosen substitution base $\gamma^k$ is sufficient to losslessly encode all of the coefficients of $r(x)$. Moreover, since $k \geq d$, the same is true of $f(x)$ and $g(x)$. Thus, the only way to have
\begin{align*}
    r(b) \not= \left(f(\gamma^k)^k \bmod{\tilde{g}(\gamma^k)}\right) \bmod{(\gamma^k - b)} ,
\end{align*}
is if
\begin{align*}
    & (\gamma^k - b) \mid (f(\gamma^k)^k \bmod{\tilde{g}(\gamma^k)}) \\
    & \implies (\gamma^k - b) \mid r(b) \\
    & \implies r(b) \equiv 0 \pmod{\gamma^k - b} .
\end{align*}
However, we are given
\begin{align*}
    r(b) \not\equiv 0 \pmod{\gamma^k - b} .
\end{align*}
\end{proof}

\begin{theorem} \label{proof:crecursiveencoding}
Let $n,d \in \Z^+$. Define $A(n)$ to be a C-recursive integer sequence of order $d$ with initial starting conditions $A(0)=A(1)=\cdots=A(d-1)=1$, such that:
\begin{align*}
    A(n) = c_{d-1} A(n-1) + c_{d-2} A(n-2) + \cdots + c_{0} A(n-d) ,
\end{align*}
where the $c_j$ are coefficients in $\Z$. Consider a polynomial
\begin{align*}
    g(x) := c_{d-1} x^{d-1} + \cdots + c_1 x + c_0 \in \Z[x] .
\end{align*}
Let $b \in \Z^+$. Suppose $b > A(n)$. Then
\begin{align*}
    A(n) = \left(b^{n-1} \bmod{(b^d - g(b))}\right) \bmod{(b-1)} .
\end{align*}
\end{theorem}
\begin{proof}
Consider the ring $R = \Z[x]/(x^d-g(x))$. By \cref{proof:recurrences}, we have that $\forall n \in \Z$, the sequence term $A(n)$ can be calculated by expanding the polynomial $f(x) := x^{n-1} \in R$ and then evaluating the result at $x=1$.

In the ring $R$, all elements are implicitly reduced modulo $(x^d - g(x))$, which enforces the relation $x^d = g(x)$. For $f(x)$, we have
\begin{align*}
    f(x) &= x^{n-1} \bmod{(x^d - g(x))} .
\end{align*}
In modular arithmetic, evaluating at $x=1$ is the same as reducing modulo $(x-1)$ (by the remainder theorem). Hence, the evaluation of $f(1) \in R$ is equal to
\begin{align*}
    f(1) &= \left(x^{n-1} \bmod{(x^d - g(x))}\right) \bmod{(x-1)} \in R .
\end{align*}

Applying \cref{proof:kroneckerqrings} by substituting with $x = b$, we arrive at the given formula
\begin{align*}
    A(n) &= \left(b^{n-1} \bmod{(b^d - g(b))}\right) \bmod{(b-1)} .
\end{align*}
By \cref{proof:kroneckerqrings}, this substitution is valid since $b > A(n)$. This completes the proof.
\end{proof}

\section{The Pell sequence} \label{section:pell}
We apply \cref{proof:kroneckerqrings} to derive a new formula for the $n$-th Pell number, which is sequence \seqnum{A000129} in the OEIS \cite{A000129}. Starting from $n=0$, the sequence $P$ begins as
\begin{align*}
    P_n = 0, 1, 2, 5, 12, 29, 70, 169, 408, 985, 2378, 5741, 13860, 33461, 80782, 195025, \ldots .
\end{align*}

\begin{theorem} \label{proof:pell}
Let $n \in \Z^+$. Let $P_n$ denote the $n$-th term of the Pell sequence, such that $P_0 = 0, P_1 = 1$, and for $n > 1$:
\begin{align*}
    P_n = 2 P_{n-1} + P_{n-2} .
\end{align*}
Then
\begin{align*}
P_n &= \left((3^n+1)^{n-1} \bmod{(9^n-2)}\right) \bmod{(3^n-1)}.
\end{align*}
\end{theorem}
\begin{proof}
Fix a ring $R = \Z[x]/(x^2 - 2)$. Consider $f(x) = x+1 \in R$. Expanding $f(x)^n$ using the binomial theorem gives
\begin{align*}
    f(x)^n = (x+1)^n = \sum_{k=0}^{n} \binom{n}{k} x^n .
\end{align*}
In the ring $R$, the elements obey the relation $x^2 = 2$. This means that in the ring $R$, all elements are implicitly reduced modulo $(x^2 - 2)$. Factoring out $x^2$ from $x^k$ in our previous expansion (to apply the modular reduction), and then simplifying, yields
\begin{align*}
    f(x)^n \bmod{(x^2 - 2)} = \sum_{k=0}^{n} \binom{n}{k} 2^{\floor{\frac{1}{2} k}} x
    = \sum_{k=0}^{n} \binom{n}{k} 2^{\floor{\frac{k}{2}}} x .
\end{align*}
Recall that evaluating at $x=1$ is the same as reducing modulo $(x-1)$. Thus, we have
\begin{align*}
    \left(f(x)^n \bmod{(x^2-2)}\right) \bmod{(x-1)} = \sum_{k=0}^{n} \binom{n}{k} 2^{\floor{\frac{k}{2}}} .
\end{align*}
From \cite{A000129}, we have the formula
\begin{align*}
    P_{n+1} = \sum_{k=0}^{n} \binom{n}{k} 2^{\floor{\frac{k}{2}}} .
\end{align*}
Hence
\begin{align*}
    P_{n+1} &= \left( f(x)^n \bmod{(x^2-2)} \right) \bmod{(x-1)} \\
    &= \left( (x+1)^n \bmod{(x^2-2)} \right) \bmod{(x-1)} .
\end{align*}
Replacing $n$ with $n-1$ yields
\begin{align*}
    P_n &= \left((x+1)^{n-1} \bmod{(x^2-2)}\right) \bmod{(x-1)} .
\end{align*}
To arrive at our stated formula, we apply \cref{proof:kroneckerqrings} by substituting with $x = 3^n$, followed by simplifying
\begin{align*}
     P_n &= \left((3^n+1)^{n-1} \bmod{((3^n)^2-2)}\right) \bmod{(3^n-1)} \\
     &= \left((3^n+1)^{n-1} \bmod{(9^n-2)}\right) \bmod{(3^n-1)} .
\end{align*}
By \cref{proof:kroneckerqrings}, this substitution is valid since $3^n > P_n$.

Considering $n = 0$, since the final modulus $(3^0-1) = 1 - 1 = 0$ results in an undefined remainder (due to division by zero). Thus, the formula is valid for $n > 0$.
\end{proof}

\section{The central binomial coefficients} \label{section:cbc}
We apply \cref{proof:kroneckerqrings} to derive a new formula for the $n$-th central binomial coefficient $\binom{2n}{n}$, which is sequence \seqnum{A000984} in the OEIS \cite{A000984}. Starting from $n=0$, the sequence of central binomial coefficients begins as
\begin{align*}
    \binom{2n}{n} = 1, 2, 6, 20, 70, 252, 924, 3432, 12870, 48620, 184756, 705432, 2704156, 10400600, \ldots
\end{align*}

\begin{theorem} \label{proof:cbc}
Let $n \in \Z^+$ such that $n > 0$. Then
\begin{align*}
\binom{2n}{n} &= \left((4^n + 1)^{2n} \bmod{(4^{n(n+1)} + 1)}\right) \bmod{(4^n-1)} .
\end{align*}
\end{theorem}
\begin{proof}
Fix a ring $R = \Z[x]/(x^{n+1} + 1)$. In the ring $R$, the elements obey the relation $x^{n+1} = -1$.

Let $f(x) = (x + 1)^{2n} \in R$. Expanding $f(x)$ and taking the result modulo $(x-1)$ gives
\begin{align*}
& \left((x + 1)^{2n} \bmod{(x^{n+1} + 1)}\right) \bmod{(x-1)} \\
&= \sum_{k=0}^{2n} \binom{2n}{k} (-1)^{\floor{\frac{k}{n+1}}} \\
&= \left( \sum_{k=0}^{n} \binom{2n}{k} (-1)^{\floor{\frac{k}{n+1}}} \right) + \left( \sum_{k=n+1}^{2n} \binom{2n}{k} (-1)^{\floor{\frac{k}{n+1}}} \right) \\
&= \left( \sum_{k=0}^{n} \binom{2n}{k} (-1)^0 \right) + \left( \sum_{k=n+1}^{2n} \binom{2n}{k} (-1)^1 \right) \\
&= \left( \sum_{k=0}^{n} \binom{2n}{k} \right) - \left( \sum_{k=n+1}^{2n} \binom{2n}{k} \right) \\
&= \binom{2n}{n} .
\end{align*}

Thus, we have
\begin{align*}
    \binom{2n}{n} &= \left((x + 1)^{2n} \bmod{(x^{n+1} + 1)}\right) \bmod{(x-1)} .
\end{align*}

Applying \cref{proof:kroneckerqrings} by substituting with $x = 4^n$ and simplifying, yields
\begin{align*}
    \binom{2n}{n} &= \left((4^n + 1)^{2n} \bmod{(4^{n(n+1)} + 1)}\right) \bmod{(4^n-1)} .
\end{align*}
By \cref{proof:kroneckerqrings}, this substitution is valid since $4^n > \binom{2n}{n}$.

Considering $n = 0$, the final modulus $(4^n-1) = 4^0-1 = 0$ results in an undefined remainder (due to division by zero). Thus, the formula is valid for $n > 0$.
\end{proof}

\section{Integer roots}
Our result on Pell numbers (\cref{proof:pell}) leads us to an interesting result on the $n$-th roots of positive integers. We require a lemma:

\begin{lemma} \label{proof:rootspolynomial}
Let $a,n \in \Z^+$. Fix a ring $S = \Z[x]/(x^n - a)$ and define the polynomial function
\begin{align*}
f_k(x) := (x+1)^k \in S .
\end{align*}
Then
\begin{align*}
\sqrt[n]{a} &= \lim_{k\to\infty} \frac{f_{k+1}(1)}{f_k(1)} - 1 \in \mathbb{R} .
\end{align*}
\end{lemma}
\begin{proof}
In $S$, the equivalence class $[x]$ satisfies the relation $[x^n] = [a]$. By the universal property of quotient rings, there exists a ring isomorphism:
\begin{align*}
S = \Z[x]/(x^n - a) \cong \Z[\sqrt[n]{a}] ,
\end{align*}
where $\Z[\sqrt[n]{a}]$ is the subring of $\mathbb{Q}(\sqrt[n]{a})$ generated by $\{ 1, \sqrt[n]{a}, \sqrt[n]{a}^2, \ldots, \sqrt[n]{a}^{n-1} \}$, and the isomorphism is induced by mapping $[x] \mapsto \sqrt[n]{a}$.

For any exponential monomial $x^k$ with $k \geq n$, we use the relation $x^k = a$ to reduce higher powers of $x$. Thus, $(x+1)^k$ can be expressed as a linear combination of $\{ 1,x,x^2,\ldots,x^d \}$ in $S$.

By the binomial theorem, we have
\begin{align*}
f_k(x) = (x+1)^k \bmod (x^n-a)
= \sum_{j=0}^k \binom{k}{j} x^j \bmod (x^n-a) 
= \sum_{j=0}^k \binom{k}{j} a^{\floor{j/n}} x^{j \bmod n} .
\end{align*}
Upon evaluation at $x=1$, this simplifies to:
\begin{align*}
f_k(1) = \sum_{j=0}^k \binom{k}{j} a^{\floor{j/n}} .
\end{align*}

The quotient
\begin{align*}
\frac{f_{k+1}(1)}{f_k(1)}
= \frac{\sum_{j=0}^{k+1} \binom{k+1}{j} a^{\floor{j/n}}}{\sum_{j=0}^k \binom{k}{j} a^{\floor{j/n}}}
\in \mathbb{R}
\end{align*}
approximates $x + 1 = \sqrt{n} + 1$ as $k \to \infty$ because the coefficients in $\frac{f_{k+1}(x)}{f_k(x)}$ become closer and closer, leading to
\begin{align*}
\left| \sqrt[n]{a} - \frac{f_{k+2}(1)}{f_{k+1}(1)} \right|
< \left| \sqrt[n]{a} - \frac{f_{k+1}(1)}{f_{k}(1)} \right| \in \mathbb{R}
\end{align*}
and an increasingly accurate approximation. Thus
\begin{align*}
\sqrt[n]{a} &= \lim_{k\to\infty} \frac{f_{k+1}(1)}{f_k(1)} - 1 \in \mathbb{R} .
\end{align*}
\end{proof}

\begin{theorem} \label{proof:roots}
Let $a,n \in \Z^+$ such that $a > 1$. Then
\begin{align*}
\sqrt[n]{a} &= \lim_{k\to\infty}
    \frac{(k^{kn} + 1)^{kn+1} \bmod{(k^{kn^2}-a)}}
    {(k^{kn} + 1)^{kn} \bmod{ (k^{kn^2}-a)}} - 1  \quad \in \mathbb{R} .
\end{align*}
\end{theorem}
\begin{proof}
By \cref{proof:rootspolynomial}, we have that
\begin{align*}
\sqrt[n]{a} &= \lim_{k\to\infty} \frac{f_{k+1}(1)}{f_k(1)} - 1 \in \mathbb{R} .
\end{align*}
As an application of \cref{proof:kroneckerqrings} with $x = k^{kn}$, we have:
\begin{align*}
\sqrt[n]{a} &= \lim_{k\to\infty}
    \frac{(k^{kn} + 1)^{kn+1} \bmod{(k^{kn^2}-a)}}
    {(k^{kn} + 1)^{kn} \bmod{ (k^{kn^2}-a)}} - 1  \quad \in \mathbb{R} .
\end{align*}
The substitution $x=k^{kn}$ is valid since $k^{kn}$ is greater than the coefficient sum of $f(x)^{kn} \bmod{(x^n-a)}$. Furthermore, the substitution does not alter the value of $x$ on the left-hand side because it is injective and only used temporarily; the value of $x$ is restored to $\sqrt[n]{a}$ upon decoding from the base $k^{kn}$, which is being carried out by the modulo operations.
\end{proof}

An intriguing aspect of the formula for $\sqrt[n]{a}$ given by \cref{proof:roots} is the applicability of a second modular reduction when evaluating the polynomial $f_k(x)$. Although it appears that both the original and modified formulas theoretically converge to the same value, our cursory analysis suggests that the convergence rates and the actual values in $\mathbb{R}$ differ. Let $c \in \Z_{\geq -1}$. Then, we observe
\begin{align*}
\sqrt[n]{a} = \lim_{k \to \infty}
\frac{\left( (k^{kn} + 1)^{kn+1} \bmod{(k^{kn^2}-a)} \right) \bmod{(k^{kn}-c)}}
    {\left( (k^{kn} + 1)^{kn} \bmod{ (k^{kn^2}-a)} \right) \bmod{(k^{kn}-c)}} - 1  \quad \in \mathbb{R},
\end{align*}
contrasted with
\begin{align*}
\sqrt[n]{a} &= \lim_{k \to \infty}
    \frac{(k^{kn} + 1)^{kn+1} \bmod{(k^{kn^2}-a)}}
    {(k^{kn} + 1)^{kn} \bmod{ (k^{kn^2}-a)}} - 1  \quad \in \mathbb{R}.
\end{align*}
It is important to note that both formulas are convergent, but not monotone due to the modular reductions. A deeper examination into the differences between these expressions and their numerical stabilities could yield further simplifications or conclusions.

From our limit formula for $\sqrt[n]{a}$, it is possible to derive a simple elementary formula for the integer part of $n$-th roots when $k$ is chosen to be sufficiently large. However, this would require a deeper examination of the convergence properties, and is out of scope for the current paper. For now, we make a conjecture:
\begin{conjecture} \label{conjecture:integerroots}
Let $a,n \in \Z^+$ such that $a > 2$, $\floor{\log_2(a)} + 1 \geq n > 1$, and $\not\exists k \in \Z^+ : k^n = a$. Then
\begin{align*}
\floor{\sqrt[n]{a}}
&= \floor{\frac{(a^{2an} + 1)^{2an+1} \bmod{(a^{2an^2}-a)}}{(a^{2an} + 1)^{2an} \bmod{ (a^{2an^2}-a)}} - 1} .
\end{align*}
\end{conjecture}


\begingroup
\raggedright
\bibliographystyle{unsrtnat}
\bibliography{main}

\begin{thebibliography}{9}
\providecommand{\natexlab}[1]{#1}
\providecommand{\url}[1]{\texttt{#1}}
\expandafter\ifx\csname urlstyle\endcsname\relax
  \providecommand{\doi}[1]{doi: #1}\else
  \providecommand{\doi}{doi: \begingroup \urlstyle{rm}\Url}\fi

\bibitem[{J. von zur Gathen and J. Gerhard}(2013)]{gathen2013modern}
{J. von zur Gathen and J. Gerhard}.
\newblock \emph{{Modern Computer Algebra}}.
\newblock Cambridge University Press, 3rd edition, 2013.
\newblock ISBN 978-1107039032.

\bibitem[{D. Harvey}(2009)]{harvey2009kronecker}
{D. Harvey}.
\newblock {Faster Polynomial Multiplication via Multipoint Kronecker
  Substitution}.
\newblock \emph{Journal of Symbolic Computation}, 44, 2009.
\newblock \doi{10.1016/j.jsc.2009.05.004}.

\bibitem[{D. Harvey and J. van der Hoeven}(2019)]{harvey2019faster}
{D. Harvey and J. van der Hoeven}.
\newblock {Faster Polynomial Multiplication Over Finite Fields Using Cyclotomic
  Coefficient Rings}.
\newblock \emph{Journal of Complexity}, 54, 2019.
\newblock ISSN 0885-064X.
\newblock URL
  \url{https://www.sciencedirect.com/science/article/pii/S0885064X19300378}.

\bibitem[{M. R. Albrecht, C. Hanser, A. Hoeller, T. Pöppelmann, F. Virdia, and
  A. Wallner}(2018)]{albrecht2018implementing}
{M. R. Albrecht, C. Hanser, A. Hoeller, T. Pöppelmann, F. Virdia, and A.
  Wallner}.
\newblock {Implementing RLWE-based Schemes Using an RSA Co-Processor}.
\newblock Cryptology ePrint Archive, Paper 2018/425, 2018.
\newblock URL \url{https://eprint.iacr.org/2018/425}.

\bibitem[{J. W. Bos, J. Renes, C. van Vredendaal}(2020)]{bos2020postquantum}
{J. W. Bos, J. Renes, C. van Vredendaal}.
\newblock {Post-Quantum Cryptography with Contemporary Co-Processors: Beyond
  Kronecker, Schönhage-Strassen and Nussbaumer}.
\newblock Cryptology ePrint Archive, Paper 2020/1303, 2020.
\newblock URL \url{https://eprint.iacr.org/2020/1303}.

\bibitem[{A. Greuet, S. Montoya, and C. Vermeersch}(2022)]{greuet2022modular}
{A. Greuet, S. Montoya, and C. Vermeersch}.
\newblock {Modular Polynomial Multiplication Using RSA/ECC coprocessor}.
\newblock Cryptology ePrint Archive, Paper 2022/879, 2022.
\newblock URL \url{https://eprint.iacr.org/2022/879}.

\bibitem[{J. M. Shunia}(2023)]{shunia2023simple}
{J. M. Shunia}.
\newblock {Simple Formulas for Univariate Multinomial Coefficients}, 2023.
\newblock URL \url{https://arxiv.org/abs/2312.00301}.
\newblock {arXiv Preprint}.

\bibitem[Inc.(2024{\natexlab{a}})]{A000129}
OEIS~Foundation Inc.
\newblock {Pell Numbers - Entry A000129 in The On-Line Encyclopedia of Integer
  Sequences}.
\newblock \url{https://oeis.org/A000129}, 2024{\natexlab{a}}.

\bibitem[Inc.(2024{\natexlab{b}})]{A000984}
OEIS~Foundation Inc.
\newblock {Central Binomial Coefficients - Entry A000984 in The On-Line
  Encyclopedia of Integer Sequences}.
\newblock \url{https://oeis.org/A000984}, 2024{\natexlab{b}}.

\end{thebibliography}
\endgroup

\end{document}